\documentclass[11pt,leqno]{amsart}

\usepackage[english]{babel}

\usepackage[letterpaper,top=2cm,bottom=2cm,left=3cm,right=3cm,marginparwidth=1.75cm]{geometry}

\usepackage{amsmath, amsthm, amssymb,latexsym, mathrsfs}
\usepackage{graphicx}
\usepackage[colorlinks=true, allcolors=blue]{hyperref}
\usepackage{xcolor}

\usepackage{tikz}
\usetikzlibrary{shapes.geometric, positioning, arrows.meta,}
\usepackage{subcaption}
\usepackage{geometry}

\newtheorem{theorem}{\bf Theorem}[section]
\newtheorem{lemma}[theorem]{\bf Lemma}
\newtheorem{example}[theorem]{\bf Example}
\newtheorem{proposition}[theorem]{\bf Proposition}
\newtheorem{corollary}[theorem]{\bf Corollary}
\newtheorem{remark}[theorem]{\sc Remark}

\newtheorem{problem}[theorem]{\bf Problem}

\usepackage{datetime}
\usepackage{float}

\DeclareMathOperator{\PSL}{PSL}

\title{Generalized Power Graph of a group}

\author{Danilo Silveira}
\address[D. Silveira]{Departamento de Ci\^encias Exatas e Aplicadas, Universidade Federal de Ouro Preto, 35931-008 	Ouro Preto, MG,  Brazil}
\address{Departamento de Matemática, Universidade Federal Minas Gerais,  31270-901	Belo Horizonte, MG,  Brazil  }
\email{danilo.sancao@ufop.edu.br}

\author{Yerko Contreras Rojas}
\address[Y. Contreras-Rojas]
{Departamento de Matemática, Universidade Estadual de Maringá 
87020-900 Maringá, PR,  Brasil }
\email{ycrojas@uem.br}

\begin{document}

\begin{abstract}
In this paper, we introduce a new class of graphs, called generalized power graphs, and investigate some of their key properties, such as completeness, connectivity, and planarity.
\end{abstract}

\maketitle
\section{Introduction}
The concepts of directed and undirected power graphs of a group were introduced by Kelarev and Quinn \cite{KQ} and Chakrabarty et al. \cite{CGS}, respectively,  and have since attracted considerable interest (see for instance \cite{AKC, AAA,  DR, DFG, KSCC, MR}). 
For the reader's convenience, we shortly define these two types of graphs. 
The \textit{directed power graph} $\vec{\mathscr{P}}(G)$ of a group $G$ is the directed graph with vertex set $G$, in which there is a directed edge from $x$ to $y$ if  $x\neq y$ and $y = x^n$ for some positive integer $n$. The undirected power graph $\mathscr{P}(G)$ of $G$ is the undirected graph obtained from $\vec{\mathscr{P}}(G)$ by ignoring edge directions. The 
proper power graph $\mathscr{P}^{*}(G)$ of $G$ is the undirected power graph without identity element of a group $G$.

An element of a group $G$ is called a generalized torsion element if a finite product of its conjugates equals the identity. We say that $G$ is a generalized torsion group if all of its elements have this property.  
It is immediate that any torsion element is a generalized torsion element. Although the converse holds for many important classes of groups (cf. \cite{BL, BSS, BS, KS}), there exist infinite groups containing non-torsion elements that are nevertheless generalized torsion. For instance, consider the non-torsion elements of the infinite dihedral group.  

Investigating generalized torsion elements provides important insights into  combinatorial group theory and topology because they serve as a primary obstruction to bi-orderability. Indeed, if a group admits a bi-invariant order, then any finite product of conjugates of a positive element remains strictly positive. Consequently, every bi-orderable group is necessarily free of generalized torsion. This connection has motivated the study of generalized torsion and has established it as an active area of research (cf.~\cite{BL, BSS, BS, I, IMT, IMT2, KS}).

Given that graphs constructed from algebraic objects often reflect the structural properties of the underlying system, following \cite{CGS, KQ}, we introduce the following classes of graphs.
The \textit{generalized directed power graph}  $\vec{\mathscr{P}}_{\bullet}(G)$ of a group $G$ is the directed graph with vertex set $G$, where for any $x, y \in G$, there is a directed edge from $x$ to $y$ (denoted $x \rightarrow y$ or $y \leftarrow x$) if $x\neq y$ and  $y$ can be expressed as a finite product of conjugates of $x$. That is, $x \rightarrow y$ if there exist $g_1, \ldots, g_k \in G$ such that$$y = x^{g_1} x^{g_2} \dots x^{g_k},$$ where $x^g = g^{-1}xg$ denotes the conjugate of $x$ by $g$. The \textit{generalized undirected power graph} (or \textit{generalized power graph}, for short)  ${\mathscr{P}}_{\bullet}(G)$ of a group $G$ is the undirected graph  obtained from $\vec{\mathscr{P}}_{\bullet}(G)$ by ignoring edge directions.
The \textit{proper generalized directed (resp., undirected) power graph} of $G$, denoted by $\vec{\mathscr{P}}^{*}_{\bullet}(G)$ (resp., ${\mathscr{P}}^{*}_{\bullet}(G)$), is the induced subgraph obtained from $\vec{\mathscr{P}}_{\bullet}(G)$ (resp., ${\mathscr{P}}_{\bullet}(G)$) by deleting the vertex corresponding to the identity element. 

 Let $S_3$ denote the symmetric group on three letters. The image below illustrates the generalized directed  power graph, the proper generalized directed power graph, and the proper generalized undirected power graph of $S_3$, respectively.

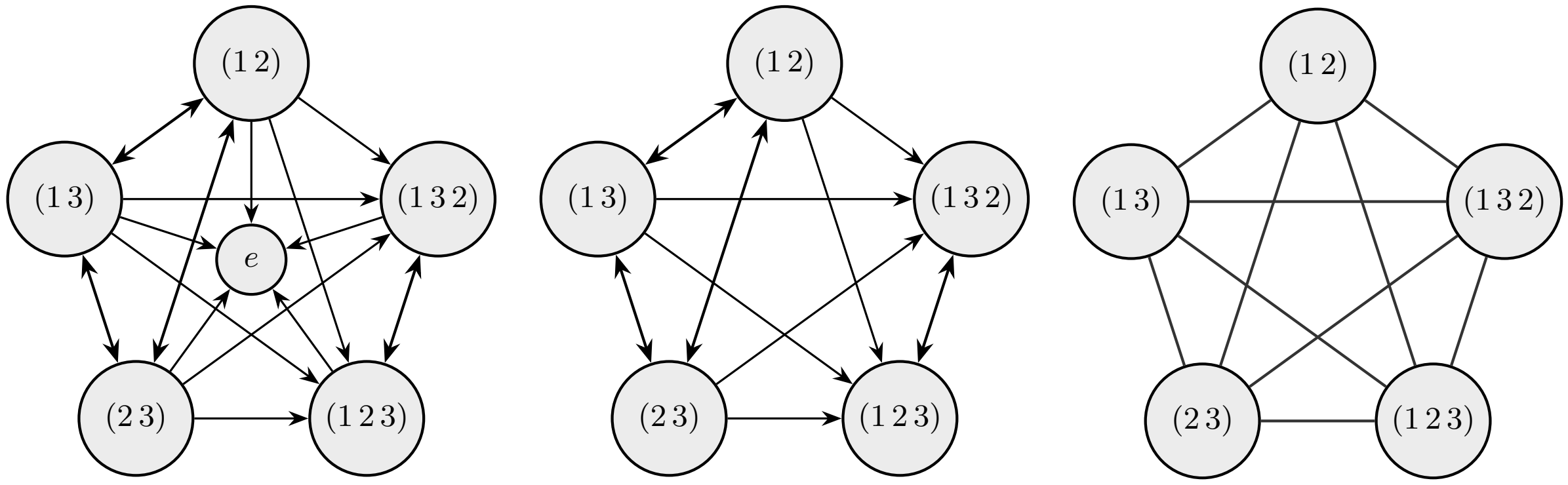
\begin{center}
\end{center}
\begin{figure}[htbp]
    \centering
    \tikzset{
        node_style/.style={circle, draw=black, fill=gray!15, inner sep=0pt, minimum size=11.5mm, font=\small, thick},
        special_node/.style={circle, draw=black, fill=gray!15, inner sep=0pt, minimum size=7mm, font=\small, thick},
        edge_undirected/.style={-, thick, black!80},
        t_to_t/.style={->, >=Stealth, black, thick},
        c_to_c/.style={->, >=Stealth, black, thick},
        t_to_c/.style={->, >=Stealth, black, semithick}
    }

    \begin{subfigure}[c]{0.31\textwidth}
        \centering
        \begin{tikzpicture}[scale=1.1]
            \node[special_node] (id)   at (0,0)       {\footnotesize$e$};
            
            \node[node_style] (t12)  at (90:1.8cm)  {\footnotesize$(1\,2)$};
            \node[node_style] (t13)  at (162:1.8cm) {\footnotesize$(1\,3)$};
            \node[node_style] (t23)  at (234:1.8cm) {\footnotesize$(2\,3)$};
            \node[node_style] (c123) at (306:1.8cm) {\footnotesize$(1\,2\,3)$};
            \node[node_style] (c132) at (18:1.8cm)  {\footnotesize$(1\,3\,2)$};

            \draw[t_to_t, <->] (t12) -- (t13);
            \draw[t_to_t, <->] (t12) -- (t23);
            \draw[t_to_t, <->] (t13) -- (t23);
            \draw[c_to_c, <->] (c123) -- (c132);

            \foreach \t in {t12, t13, t23} {
                \draw[t_to_c] (\t) -- (c123);
                \draw[t_to_c] (\t) -- (c132);
            }

            \foreach \v in {t12, t13, t23, c123, c132} {
                \draw[t_to_c] (\v) -- (id);
            }
        \end{tikzpicture}
         \caption{$\vec{\mathscr{P}}_{\bullet}(S_3)$}
    \end{subfigure}\hfill
    %
    \begin{subfigure}[c]{0.31\textwidth}
        \centering
        \begin{tikzpicture}[scale=1.1]
            \node[node_style] (t12)  at (90:1.8cm)  {\footnotesize$(1\,2)$};
            \node[node_style] (t13)  at (162:1.8cm) {\footnotesize$(1\,3)$};
            \node[node_style] (t23)  at (234:1.8cm) {\footnotesize$(2\,3)$};
            \node[node_style] (c123) at (306:1.8cm) {\footnotesize$(1\,2\,3)$};
            \node[node_style] (c132) at (18:1.8cm)  {\footnotesize$(1\,3\,2)$};

            \draw[t_to_t, <->] (t12) -- (t13);
            \draw[t_to_t, <->] (t12) -- (t23);
            \draw[t_to_t, <->] (t13) -- (t23);
            \draw[c_to_c, <->] (c123) -- (c132);

            \foreach \t in {t12, t13, t23} {
                \draw[t_to_c] (\t) -- (c123);
                \draw[t_to_c] (\t) -- (c132);
            }
        \end{tikzpicture}
        \caption{$\vec{\mathscr{P}}^{\ast}_{\bullet}(S_3)$}
    \end{subfigure}\hfill
    %
    \begin{subfigure}[c]{0.31\textwidth}
        \centering
        \begin{tikzpicture}[scale=1.1]
            \node[node_style] (t12)  at (90:1.8cm)  {\footnotesize$(1\,2)$};
            \node[node_style] (t13)  at (162:1.8cm) {\footnotesize$(1\,3)$};
            \node[node_style] (t23)  at (234:1.8cm) {\footnotesize$(2\,3)$};
            \node[node_style] (c123) at (306:1.8cm) {\footnotesize$(1\,2\,3)$};
            \node[node_style] (c132) at (18:1.8cm)  {\footnotesize$(1\,3\,2)$};

            \draw[edge_undirected] (t12) -- (t13);
            \draw[edge_undirected] (t12) -- (t23);
            \draw[edge_undirected] (t13) -- (t23);
            \draw[edge_undirected] (c123) -- (c132);
            \foreach \t in {t12, t13, t23} {
                \draw[edge_undirected] (\t) -- (c123);
                \draw[edge_undirected] (\t) -- (c132);
            }
        \end{tikzpicture}
        \caption{$\mathscr{P}^{\ast}_{\bullet}(S_3)$}
    \end{subfigure}

\end{figure}

Given that every non-trivial element is adjacent to the identity in the generalized power graph of a generalized torsion group, our results are mainly devoted to the proper generalized power graph.

The paper is organized as follows. In Section 2, we present general results regarding generalized power graphs. Sections 3 and 4 are devoted to investigating the completeness and connectedness of the generalized power graph of some families of finite and infinite groups. Finally, in Section 5, we give necessary and sufficient conditions for the planarity of the proper generalized power graph of a finite group.

\section{General Properties}

Let $x$ and $y$ be elements in a group $G$. We say that there exists a path connecting $x$ and $y$ if there exist elements $x=z_1, z_2, \ldots, z_n= y \in G$ such that $z_i$ and $z_{i+1}$ are adjacent for all $i = 1, \ldots, n-1$.
 
For non-empty subsets $X, Y \subseteq G$, we write $X \to Y$ if $x \to y$ for all $x \in X$ and $y \in Y$. We may also write $X \leftrightarrow Y$ to mean that $X \to Y$ and $Y \to X$.

We denote the conjugacy class of $x$ in $G$ by $x^G = \{x^g \mid g \in G\}$, $|x|$ the order of $x$ and its centralizer by $C_G(x)$. Furthermore, $C_n, D_{2n}$, $Q_8$, $K_n$, and $K_{3,3}$ denote the cyclic group of order $n$, the dihedral group of order $2n$, the quaternion group of order 8, the complete graph on $n$ vertices, and the complete bipartite graph, respectively.

We define commutators of weight $1, 2, \dots$ in the variables $x_1, x_2, \dots$ recursively as formal bracket expressions. The variables $x_1, x_2, \dots$ are commutators of weight $1$, and commutators of weight $2$ take the form $[x_i, x_j] = x_i^{-1}x_j^{-1}x_i x_j$. If $w_1$ and $w_2$ are commutators of weight $k_1$ and $k_2$, respectively, then $[w_1, w_2]$ is a commutator of weight $k_1 + k_2$.
If $w$ is a commutator of weight $k$, we view it primarily as a function of $k$ variables defined on a group $G$. Any value of $w$ in $G$ is called a commutator in $G$ (of weight $k$).

Throughout this paper, the following properties may be used without explicit reference.

\begin{lemma}\label{L1}  Let $x,y,z$ elements of a group $ G$.
\begin{enumerate}
    \item    If $x\to y \rightarrow z$, then $x\to z$. 
    \item If $x\to y$ and $x\to z$, then $x\to yz.$ 
    \item $x^G\leftrightarrow x^G$. 
    \item If $x^g\to y^h$ for some $g,h\in G$, then $x^G\rightarrow y^G$.

    If additionally $x$ is a torsion element, we obtain the following results:
     \item If $y\neq 1$ and $y\in \langle x^G \rangle$, then $x^G \rightarrow  y^G $. 
      \item If $(|x|,i)=1$, then $x^i\to x^j$ for every $j.$
    \item   $x\leftrightarrow x^{-1}$ and $y\to (xy)^{|x|}$. 
        \item If  $x$ is any entry of some commutator $w$, then $x\to w$.
\end{enumerate}
\end{lemma}

\begin{proof}
\begin{enumerate}
      \item  There exist elements $g_1,\ldots,g_m, h_1,\ldots,h_n \in G$ such that $y=x^{g_1}\cdots x^{g_m}$ and $z=y^{h_1}\cdots y^{h_n}$. Hence,$$z = \prod_{j=1}^n (x^{g_1}\cdots x^{g_m})^{h_j} = \prod_{j=1}^n \prod_{i=1}^m x^{g_i h_j},$$ showing that $x \to z$.
    
    \item There exist elements $g_1,\ldots,g_m, h_1,\ldots,h_n \in G$ such that $y=x^{g_1}\cdots x^{g_m}$ and $z=x^{h_1}\cdots x^{h_n}$. Hence, $yz = x^{g_1}\cdots x^{g_m} x^{h_1}\cdots x^{h_n}$, showing that $x \to yz$.

     \item Let $z_1, z_2 \in x^G$. There exist $g, h \in G$ such that $z_1 = x^g$ and $z_2 = x^h$. Thus, $z_1 = (x^h)^{h^{-1}g} = z_2^{h^{-1}g}$, which means $z_2 \to z_1$. By a similar argument, we obtain $z_1 \to z_2$.
    
       \item Note that $x^g \to y^h$ implies $x \to y$. Applying Item (3), we obtain $x^G \to \{x\} \to \{y\} \to y^G$, and the result follows from Item (1).

      \item    Since $x$ is torsion and $y \in \langle x^G \rangle$, we have $x \to y$. Thus, the result follows by Item (4).

      \item By Bézout's identity, there exist integers $m$ and $n$ such that $mi + n\vert{}x\vert{} = 1$. Thus, $x = x^{mi + n\vert{}x\vert{}} = x^{mi}$, and so $x^j = (x^{mi})^j = (x^i)^{mj}$. This implies that $x^i \to x^j$.
    
      \item If $x$ has finite order $n$, then $x^{-1}$ also has order $n$ and, consequently, $x^{n-1} = x^{-1}$ and $(x^{-1})^{n-1} = x$.
    
  The second part follows from the fact that, for any integer $n \geq 2$, we have $(xy)^n = x^n y^{x^{n-1}} y^{x^{n-2}} \cdots y^x y$. 
Setting $|x| = n$, we have $(xy)^n = y^{x^{n-1}} y^{x^{n-2}} \cdots y^x y$, which means $y \to (xy)^n$.

      \item We proceed by induction on the length of the commutator $w$. Let $k$ be the order of $x$. If $w$ has length $2$, then $w = [x,y]$ or $w = [y,x]$. In the first case, we have $w = x^{-1}y^{-1}xy = x^{k-1}x^y$, and thus $x \to w$. The second case is similar. Now assume the result holds for any commutator of length less than $t$, and let $w$ be a commutator of length $t$. Then $w = [w_1,w_2]$, where $w_1$ and $w_2$ are commutators of lengths $t_1$ and $t_2$, respectively, satisfying $t_1 + t_2 = t$. Furthermore, $x$ appears as an entry in at least one of the words $w_1, w_2$. Without loss of generality, assume that $x$ is an entry of $w_1$. By the induction hypothesis, we have $x \to w_1$, which means $w_1 = x^{g_1}\cdots x^{g_n}$ for some $g_1,\ldots,g_n \in G$. Thus, $$w = [w_1,w_2] = w_1^{-1}w_1^{w_2} =  ((x^{k-1})^{g_n}\cdots (x^{k-1})^{g_1})(x^{g_1w_2}\cdots x^{g_nw_2}).$$ This shows that $x \to w$.
\end{enumerate}
\end{proof}

We note that if ${\mathscr{P}}_{\bullet}(G)$ is complete, then $G$ is a generalized torsion group. Furthermore, if $G$ is a generalized torsion group, then ${\mathscr{P}}_{\bullet}(G)$ is complete if and only if ${\mathscr{P}}^*_{\bullet}(G)$ is complete. Obviously, we cannot replace ``complete'' with ``connected'' (as shown by $G=C_2\times C_2$).

\begin{lemma}\label{lem:properties-compl} Let $G$ be a group.
    \begin{enumerate}
      \item If ${\mathscr{P}}^{*}_{\bullet}(G)$ is complete (resp. connected), then ${\mathscr{P}}^{*}_{\bullet}(H)$ is complete (resp. connected) for every direct factor $H$ of $G$.
      
        \item If ${\mathscr{P}}^{*}_{\bullet}(G)$ is complete, then ${\mathscr{P}}^{*}_{\bullet}(G/N)$ is complete for every normal subgroup $N$ of $G$.
    \end{enumerate}
\end{lemma}
\begin{proof}
\begin{enumerate}
       \item  
    Let $H$ be a direct factor of $G$. Thus, there exists a subgroup $K \le G$ such that $G = H \times K$. Identify $H$ with the subgroup $H \times \{1\} \le G$, and choose $(h_1, 1), (h_2, 1) \in H \times \{1\}$ non-trivial elements. 
    Since $\mathscr{P}^*_{\bullet}(G)$ is complete, we may assume without loss of generality that $(h_2,1) \to (h_1,1)$. Thus, there exist elements $(x_1, y_1), \ldots, (x_n, y_n) \in G$ such that $$(h_1, 1) = (h_2, 1)^{(x_1, y_1)} \cdots (h_2, 1)^{(x_n, y_n)} = \left( h_2^{x_1} \cdots h_2^{x_n}, 1 \right),$$ which implies that $h_1 = h_2^{x_1} \cdots h_2^{x_n}$. Consequently, $\mathscr{P}^{*}_{\bullet}(H)$ is complete.

  The connectedness part follows directly from the fact that if $G = H \times K$ and $(h_1, k_1)$ is adjacent to $(h_2, k_2)$, then either $h_1 \rightarrow h_2$ and $k_1 \rightarrow k_2$, or $h_2 \to h_1$ and $k_2 \rightarrow k_1$.
    
    \item Let $xN, yN$ be non-trivial distinct elements of $G/N$. Since ${\mathscr{P}}^{*}_{\bullet}(G)$ is complete, we have $x\rightarrow y$ or $y\rightarrow x$. Assume, without loss of generality, that $x \rightarrow y$. Then there exist elements $g_1, \dots, g_k \in G$ such that $y = x^{g_1} \cdots x^{g_k}.$ Projecting this equality onto the quotient $G/N$, we obtain
$yN = (xN)^{g_1N} \cdots (xN)^{g_kN},$ which completes the proof.
\end{enumerate} \end{proof}

Note that in Lemma \ref{lem:properties-compl}~(2), we cannot replace ``complete'' with ``connected''. For instance, the proper generalized power graph of the quaternion group $Q_8$ is connected, but $Q_8$ has a quotient isomorphic to the Klein four-group, whose proper generalized power graph is disconnected.

\begin{remark} \label{remark:order8}  As shown in \cite{Cam}, any two finite abelian groups with isomorphic proper power graphs must be isomorphic. 
However ${\mathscr{P}}^{*}_{\bullet}(G) \simeq {\mathscr{P}}^{*}_{\bullet}(H)$ does not imply $G \simeq H$, even when $G$ and $H$ are nilpotent.  Indeed, consider the quaternion group $Q_8$ and the dihedral group $D_8$.

\begin{figure}[htbp]
\centering
\begin{subfigure}[b]{0.48\textwidth}
    \centering
    \begin{tikzpicture}[scale=1.1, 
        every node/.style={circle, draw, fill=gray!15, inner sep=0pt, minimum size=8mm, font=\small},
        center node/.style={circle, draw, fill=gray!15, inner sep=0pt, minimum size=8mm, font=\small},
        edge style/.style={-, thick}]

        \node[center node] (m1) at (0,0) {$-1$};

        \node (i)  at (90:1.5) {$i$};
        \node (mi) at (30:1.5) {$-i$};

        \node (j)  at (150:1.5) {$j$};
        \node (mj) at (210:1.5) {$-j$};

        \node (mk) at (270:1.5) {$-k$};
        \node (k)  at (330:1.5) {$k$};

        \draw[edge style] (m1) -- (i);
        \draw[edge style] (m1) -- (mi);
        \draw[edge style] (i)  -- (mi);

        \draw[edge style] (m1) -- (j);
        \draw[edge style] (m1) -- (mj);
        \draw[edge style] (j)  -- (mj);

        \draw[edge style] (m1) -- (k);
        \draw[edge style] (m1) -- (mk);
        \draw[edge style] (k)  -- (mk);
    \end{tikzpicture}
    \caption{$\mathscr{P}^{\ast}_{\bullet}(Q_8)$}
    \label{fig:q8}
\end{subfigure}
\hfill
\begin{subfigure}[b]{0.48\textwidth}
    \centering
    \begin{tikzpicture}[scale=1.1, 
        every node/.style={circle, draw, fill=gray!15, inner sep=0pt, minimum size=8mm, font=\small},
        center node/.style={circle, draw, fill=gray!15, inner sep=0pt, minimum size=8mm, font=\small},
        edge style/.style={-, thick, >=stealth}]

        \node[center node] (r2) at (0,0) {$a^2$};

        \node (r)  at (90:1.5) {$a$};
        \node (r3) at (30:1.5) {$a^3$};

        \node (s)   at (150:1.5) {$b$};
        \node (sr2) at (210:1.5) {$a^2b$};

        \node (sr3) at (270:1.5) {$a^3b$};
        \node (sr)  at (330:1.5) {$ab$};

        \draw[edge style] (r2) -- (r);
        \draw[edge style] (r2) -- (r3);
        \draw[edge style] (r)  -- (r3);

        \draw[edge style] (r2) -- (s);
        \draw[edge style] (r2) -- (sr2);
        \draw[edge style] (s)  -- (sr2);

        \draw[edge style] (r2) -- (sr);
        \draw[edge style] (r2) -- (sr3);
        \draw[edge style] (sr) -- (sr3);
    \end{tikzpicture}
    \caption{$\mathscr{P}^{\ast}_{\bullet}(D_8)$}
    \label{fig:d8}
\end{subfigure}
\caption{Proper Generalized Power Graphs of $Q_8$ and $D_8$.}
\label{fig:power_graphs}
\end{figure}
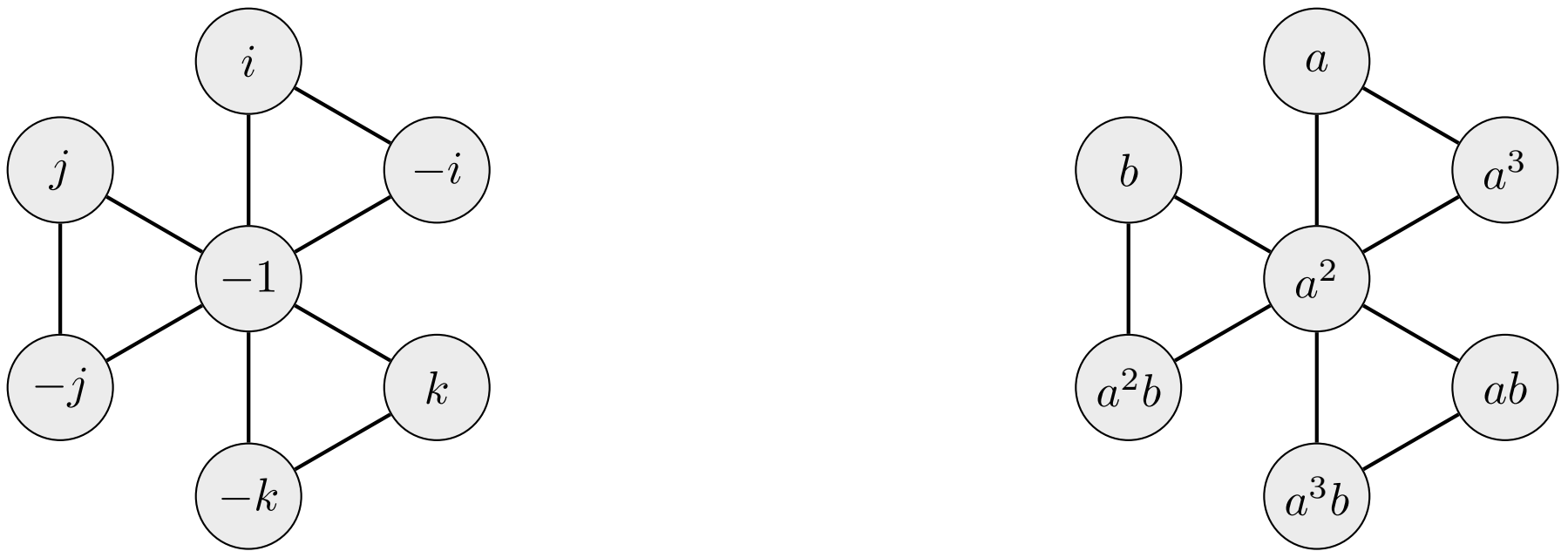

Further, the same holds if we replace ``nilpotent'' with ``finite simple'', because the proper generalized power graph of a finite simple group is complete and there exist pairs of non-isomorphic finite simple groups of the same order, such as $A_8$ and $\PSL_3(4)$.
\end{remark}

\section{Completeness}

In this section, we present several results concerning the completeness of the (proper) generalized power graph. We begin with the following result, which characterizes all abelian groups whose generalized power graph is complete.

Throughout this paper, the Prüfer $p$-group is denoted by $C_{p^{\infty}}$, where $p$ is a prime.

\begin{lemma}\label{lema:comprime} 
Let $G$ be an abelian group.  Then ${\mathscr{P}}_{\bullet}(G)$ is complete if only if $G\simeq C_{p^n}$, where $n\in \mathbb{N}\cup \{\infty\}$.
\end{lemma}
\begin{proof}
Since ${{\mathscr{P}}_{\bullet}}(G)$ coincides to power graph of $G$, the result follow by \cite[Theorem 2.12]{CGS} and \cite[Theorem 1]{AAA}.
\end{proof}

In the non-abelian case, the situation becomes much more interesting. For example, $\mathscr{P}_{\bullet}(G)$ is complete if $G=S_3$ or $A_4$ (cf. Section 1 and Example \ref{ex:A_4}), whereas $\mathscr{P}_{\bullet}(G)$ is not complete for any nonabelian group $G$ of order 8 (see Figures~\ref{fig:q8} and \ref{fig:d8} in Remark \ref{remark:order8}). We begin by proving the following result.

\begin{lemma}
If $G$ is a finitely generated nilpotent group such that $\mathscr{P}^*_\bullet(G)$ is complete, then $G$ must be a finite $p$-group.
 \end{lemma}
\begin{proof}
First, we show that $G$ must be a finite group. We proceed by induction on the nilpotency class of $G$, denoted by $c$. If $c = 1$, the result follows directly from Lemma \ref{lema:comprime}. Assume $c \ge 2$. Combining Lemma \ref{lem:properties-compl} (2) with the induction hypothesis, we deduce that the central quotient $G/Z(G)$ is a finite  group. By Schur's Theorem \cite[10.1.4]{Rob}, the derived subgroup $G'$ is finite. Lemmas \ref{lem:properties-compl} (2) and \ref{lema:comprime} ensure that the abelian quotient $G/G'$ is a finite $p$-group, which implies that $G$ is finite. Thus, $G$ decomposes as a direct product of its Sylow $p_i$-subgroups, say $G \cong P_1 \times \cdots \times P_n$. Suppose for a contradiction that $n \ge 2$ and that $\mathscr{P}^*_{\bullet}(G)$ is complete. Consider the non-trivial elements $x = (a, 1, \dots, 1)$ and $y = (1, b, \dots, 1)$ in $G$. Since $\mathscr{P}^*_{\bullet}(G)$ is complete, $x$ and $y$ must be adjacent,  say $x \to y$. Then there exist elements $g_1, \dots, g_k \in G$ such that $y = x^{g_1} \cdots x^{g_k}$. Consequently $b=1$, which is a contradiction. Therefore, $n = 1$, and $G$ is a $p$-group for some prime $p$.
\end{proof}

\begin{example}
Let $D_{2n} = \langle a, b \mid a^n = 1, b^2 = 1, bab = a^{-1} \rangle$ be the dihedral group of order $2n$. We show that $\mathscr{P}^*_{\bullet}(D_{2n})$ is complete if and only if $n = p^m$ for an odd prime $p$ and a positive integer $m$. Indeed, for every element $a^i \in \langle a \rangle$, we have $\langle a^i \rangle = \langle (a^i)^{D_{2n}} \rangle$. Consequently, if $\mathscr{P}^*_{\bullet}(D_{2n})$ is complete, then the proper power graph $\mathscr{P}^*(\langle a \rangle)$ is also complete. According to \cite[Theorem 2.12]{CGS}, this occurs only when $n = p^m$ for some prime $p$ and some positive integer $m$.

If $p$ is an odd prime, all elements $a^i b \in D_{2p^m}$ belong to the same conjugacy class $b^{D_{2p^m}}$, and $\langle b^{D_{2p^m}} \rangle = D_{2p^m}$, which implies the completeness of $\mathscr{P}^*_{\bullet}(D_{2p^m}).$

If $p=2$, it suffices to observe that the vertices $a$ and $a^2 b$ are not adjacent. Indeed, this follows from the fact that $\langle a^{D_{2\cdot 2^m}}\rangle = \langle a \rangle$ and $\langle (a^2 b)^{D_{2\cdot 2^m}}\rangle = \{a^{2k} b^i \mid 0 \le k < 2^{m-1} \text{ and } i=0,1\}$.
\end{example}

\begin{example} \label{ex:A_4} The graph ${{\mathscr{P}}^{*}_{\bullet}}(A_4)$ is complete. Indeed, the four conjugacy classes of the alternating group $A_4$ are: $C_1=\{(1)\},$
$C_2= \{(12)(34),(13)(24),(14)(23)\},$
$C_3=\{(123),(243),(134),(142)\}$ and
 $C_4=\{(132),(234),(143),(124)\}$. 
 Observe that 
     $(132)(132)^{(124)}=(142)$ and $(123)(123)^{(234)}(123)=(12)(34)$. 
Thus, Lemma \ref{L1} (4) implies that $C_4 \rightarrow C_3 \rightarrow C_2$. Combining items (1) and (3) from Lemma \ref{L1}, we obtain that ${{\mathscr{P}}^{*}_{\bullet}}(A_4)$ is complete. 
\end{example}


\begin{figure}[htbp]
\centering
\begin{tikzpicture}[scale=2, 
    v_double/.style={circle, draw, fill=gray!15, inner sep=0pt, minimum size=10mm, font=\tiny},
    v_cycle/.style={circle, draw, fill=gray!15, inner sep=0pt, minimum size=10mm, font=\tiny}]

    \node[v_double] (d1) at ({360/11 * 0}:1.8) {$(12)(34)$};
    \node[v_double] (d2) at ({360/11 * 1}:1.8) {$(13)(24)$};
    \node[v_double] (d3) at ({360/11 * 2}:1.8) {$(14)(23)$};

    \node[v_cycle] (t1) at ({360/11 * 3}:1.8) {$(123)$};
    \node[v_cycle] (t2) at ({360/11 * 4}:1.8) {$(132)$};
    \node[v_cycle] (t3) at ({360/11 * 5}:1.8) {$(124)$};
    \node[v_cycle] (t4) at ({360/11 * 6}:1.8) {$(142)$};
    \node[v_cycle] (t5) at ({360/11 * 7}:1.8) {$(134)$};
    \node[v_cycle] (t6) at ({360/11 * 8}:1.8) {$(143)$};
    \node[v_cycle] (t7) at ({360/11 * 9}:1.8) {$(234)$};
    \node[v_cycle] (t8) at ({360/11 * 10}:1.8) {$(243)$};

    \foreach \i in {d1,d2,d3,t1,t2,t3,t4,t5,t6,t7,t8} {
        \foreach \j in {d1,d2,d3,t1,t2,t3,t4,t5,t6,t7,t8} {
            \ifx\i\j\else
                \draw[gray!60, thin] (\i) -- (\j);
            \fi
        }
    }

\end{tikzpicture}
\caption{Proper Generalized Power Graph of $A_4$}
\end{figure}
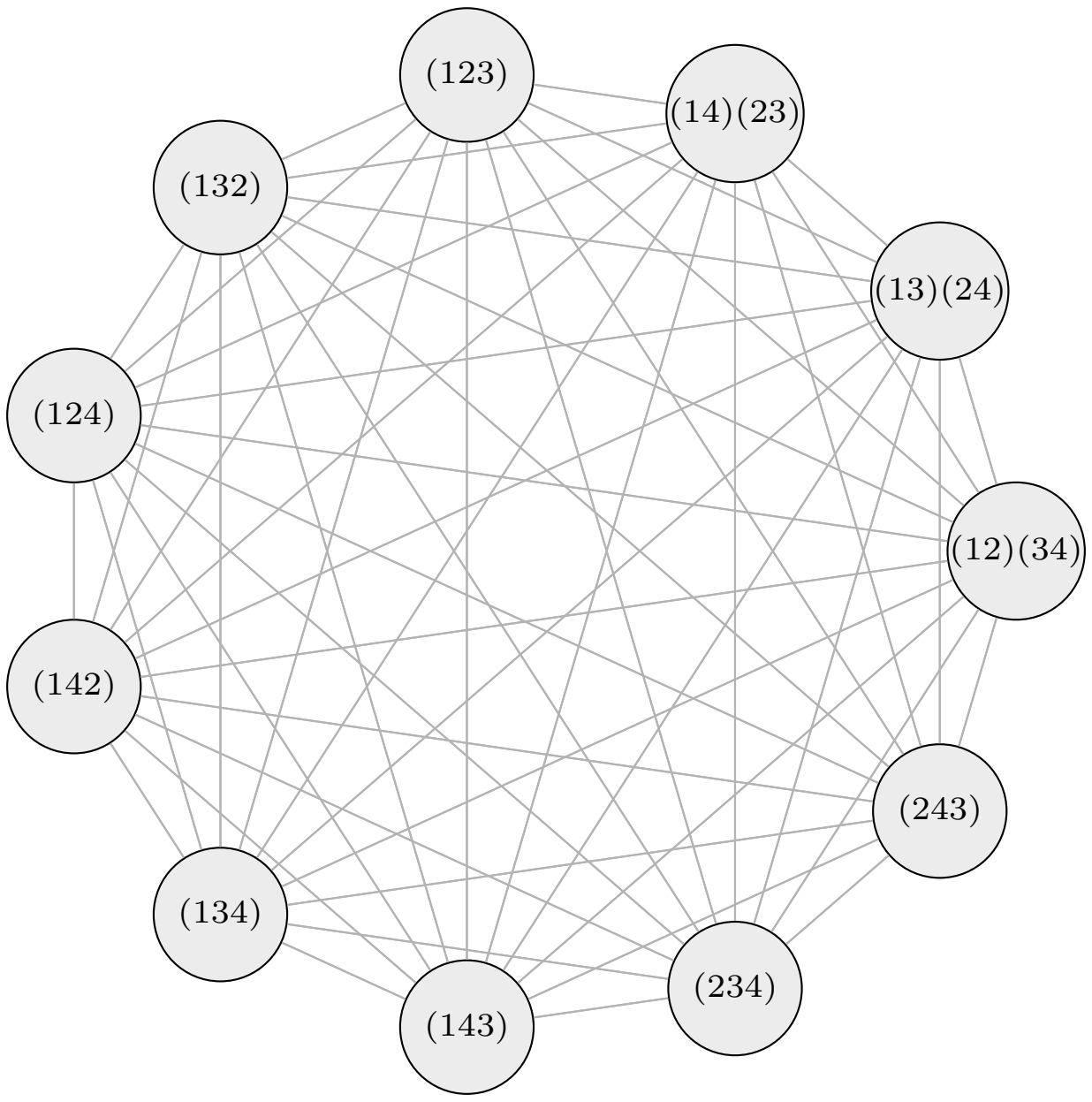

\begin{lemma}\label{lem:tsimple}
    If $G$ is a torsion simple group, then ${\mathscr{P}}^{*}_{\bullet}(G)$ is complete.
\end{lemma}
\begin{proof}
    If $G$ is abelian, we can deduce that 
$G$ has prime order. Thus, the result follows from Lemma~\ref{lema:comprime}. 
Now, assume that $G$ is nonabelian. Let $x$ and $y$ be arbitrary non-trivial elements of $G$. Since $G$ is simple, we have $G = \langle x^G \rangle$, which implies that $y = x^{\pm g_1} \cdots x^{\pm g_k}$ for some $g_1, \dots, g_k \in G$. Since $G$ is a torsion group, $x^{-1} = x^n$ for some positive integer $n$, which shows that $x \to y$. This completes the proof.
\end{proof}

\begin{remark} Completeness is hereditary under direct summands and quotients, as shown in Lemma \ref{lem:properties-compl}. However, it is not hereditary under taking subgroups. Indeed, while Lemma \ref{lem:tsimple} ensures that ${{\mathscr{P}}^*_{\bullet}}(A_5)$ is complete, the alternating group $A_5$ contains the Klein four-group $V_4 \cong C_2 \times C_2$ as a subgroup, whose corresponding graph ${{\mathscr{P}}_{\bullet}}^*(V_4)$ is not complete by Lemma \ref{lema:comprime}.
\end{remark}

As a consequence of Lemma \ref{lem:tsimple}, we obtain the following result.

\begin{theorem}\label{thm:S_n-complete} 
Every finite group can be embedded into a finite group $G$ such that ${\mathscr{P}}^{*}_{\bullet}(G)$ is complete.
   
\end{theorem}

\begin{proof}
It is sufficient to show that $\mathscr{P}^*_{\bullet}(S_n)$ is complete for all integers $n \ge 5$.
 Let $x, y \in S_n$. We prove that $y \to x$ or $x \to y$ by analysing the following cases:

\textbf{Case 1:} $x, y \in A_n$

Since $A_n$ is a nonabelian simple group for $n \geq 5$, Lemma \ref{lem:tsimple} applies.

\textbf{Case 2:} $x \in A_n$ and $y \notin A_n$

First, suppose that $y^2 \neq 1$. Observe that $y^2 \in A_n$, and by Lemma \ref{lem:tsimple}, $y^2 \leftrightarrow x$. Hence, there exist $g_1, \dots, g_m \in S_n$ such that
$$
x = (y^2)^{g_1} \cdots (y^2)^{g_m} = (y^{g_1})^2 \cdots (y^{g_m})^2,
$$
which shows that $y \to x$.

Now suppose that $y^2 = 1$. If $A_n \le C_{S_n}(y)$, since $A_n$ has index 2 in $S_n$, we get that either $C_{S_n}(y) = S_n$, which implies $y \in Z(S_n) = \{1\}$, or $C_{S_n}(y) = A_n$, which implies $y \in A_n$. Both cases yield a contradiction. Thus, there exists an element $z \in A_n$ such that $1 \neq [y,z] = y y^z \in A_n$.This implies that $y \to [y,z]$. By Lemma \ref{lem:tsimple}, we have $y \to [y,z] \leftrightarrow x$. The result follows by applying Lemma \ref{L1}(1).

\textbf{Case 3:} $y \in A_n$ and $x \notin A_n$

Repeating the argument used in the proof of Case 2, we obtain $x \to y$.

\textbf{Case 4:} $x, y \notin A_n$

Since $xy \in A_n$, by {Case 2}, we have $y \to xy$. Thus, there exist $h_1, \dots, h_s \in S_n$ such that
$$
xy = y^{h_1} \cdots y^{h_s}.
$$
Let $r$ be the order of $y$. Multiplying the right-hand side of the equality by $y^{r-1}$, we get
$$
x = y^{h_1} \cdots y^{h_s} y^{r-1},
$$
which proves that $y \to x$. 
This completes the proof.
\end{proof}

\begin{theorem}
 Let $F$ be a free group on countably many generators. Then  ${\mathscr{P}}^{\ast}_{\bullet}(F)$ is not complete. 
\end{theorem}
\begin{proof}
By Lemma \ref{lema:comprime}  we can assume that the rank of $G$ is at least 2. Assume by contradiction that ${\mathscr{P}}^{\ast}_{\bullet}(F)$ is complete. By Lemma \ref{lem:properties-compl} (2), ${\mathscr{P}}^{\ast}_{\bullet}(F/F')$ is complete as well. However, $F/F'$ is isomorphic to Cartesian product of copies of $\mathbb{Z}$ what is a contradiction (by Lemma \ref{lema:comprime}).
\end{proof}

\begin{corollary}
If $G$ has any factor direct or quotient which is isomorphic to free group on countably many generators, then ${\mathscr{P}}^{*}_{\bullet}(G)$ is not complete.
\end{corollary}
\begin{proof}
    This follows immediately from the above theorem together with Lemma \ref{lem:properties-compl}.
\end{proof}

The next lemma is probably well-known but we do not have a reference to it.

\begin{lemma}\label{lem:metacyclic}
    A finite group $G$ satisfies exactly one of the following conditions:
\begin{itemize}
    \item $G$ is metacyclic; or
    \item $G$ is non-metacyclic and contains a subgroup isomorphic to either $C_p \times C_p$ (for some prime $p$) or quaternion group of order eight $Q_8$.
\end{itemize}
\end{lemma}
\begin{proof}
Assume that $G$ is non-metacyclic and that $G$ contains no subgroup isomorphic to $C_p \times C_p$ for any prime $p$.
Since $G$ is non-metacyclic, it follows from \cite[Theorem 7.6.2]{Gor}  that some Sylow $p$-subgroup of $G$ is non-cyclic. Let $P$ be such a subgroup. This $P$ possesses a central subgroup of order $p$. Furthermore, this central subgroup is necessarily unique, otherwise $G$ would contain a copy of $C_p \times C_p$. Hence, by Burnside's theorem (cf. \cite[5.3.6]{Rob}), $P$ admits a subgroup isomorphic to the quaternion group $Q_8$.
\end{proof}

Recall that a Dedekind group is a group in which all subgroups are normal. We show the following result.

\begin{proposition}\label{lem:properties}
    \begin{enumerate}
        \item If $G$ is a Dedekind group, then ${\mathscr{P}}_{\bullet}(G)$ is complete if only if $G\simeq C_{p^n}$, where $n\in \mathbb{N}\cup \{\infty\}$.
        
        \item If $G$ is a finite group which ${\mathscr{P}}_{\bullet}(H)$ is complete for every proper subgroup and quotient $H$ of $G$, then:
        \begin{enumerate}
            \item $G$ is cyclic of order $pq$ for some different primes $p$ and $q$;
            \item $G$ is cyclic of order power of a prime;
            \item $G$ is metacyclic and is an extension of $C_{p^a}$ by $C_{q^b}$, for some primes $p, q$ and positive integers $a, b$.
        \end{enumerate}
    \end{enumerate}
\end{proposition}

\begin{proof}
\begin{enumerate}
    \item By \cite[5.3.7]{Rob}, a group $G$ is Dedekind if and only if it is either abelian or isomorphic to the direct product $Q_8 \times A$, where $A$ is a torsion abelian group. If $G$ is abelian, then the result follows by Lemma \ref{lema:comprime}. Now assume that $G = Q_8 \times A$. It follows from Lemma \ref{lem:properties-compl} (1) that if ${\mathscr{P}}_{\bullet}(G)$ is complete, then ${{\mathscr{P}}_{\bullet}}(Q_8)$ must also be complete. However,  ${{\mathscr{P}}_{\bullet}}(Q_8)$ is not complete. Therefore, in this case, ${\mathscr{P}}_{\bullet}(G)$ cannot be complete.

    \item If $G$ is non-metacyclic, then by Lemma \ref{lem:metacyclic}, $G$ contains a subgroup isomorphic to either $C_p \times C_p$ or $Q_8$. Thus, combining Item (1),  Lemma \ref{lema:comprime}  and Lemma \ref{lem:properties-compl} (2), we get a contradiction.  Now,  assume that $G$ is metacyclic. If $G$ is cyclic, then $G \simeq C_n$, where $n = p_1^{n_1} \cdots p_r^{n_r}$ and $p_1, \ldots, p_r$ are pairwise distinct primes. It follows from Lemma \ref{lema:comprime}  that $G$ is either of prime-power order or of order $pq$ for distinct primes $p$ and $q$. Assuming that $G$ is not cyclic, there exists $N \triangleleft G$ be such that both $N$ and $G/N$ are non-trivial cyclic. By Lemma \ref{lema:comprime}, $N$ and $G/N$ have prime-power orders, say $p^a$ and $q^b$, respectively. This completes the proof.
\end{enumerate}
\end{proof}

The structure of the $p$-subgroups of a finite group can reveal significant 
information about the group itself. For instance, Lucchini proved in~\cite{Lucc} 
that if every $p$-subgroup of a group $G$ can be generated by at most $r$ 
elements for every prime $p$ dividing the order of $G$, then every subgroup of $G$ can be generated by at most $r+1$ elements.
This insight naturally motivates the following problem.

\begin{problem}\label{Prob:Lucc-type}
  Let $G$ be a finite group and let $\mathcal{P}$ be a graph property. If 
${{\mathscr{P}}_{\bullet}}(H) $ satisfies $\mathcal{P}$ for every $p$-subgroup $H$ of $G$ and  every prime $p$ dividing the order of $G$, what can be said about the structure of $G$?
\end{problem}

We establish the following result.

\begin{proposition}\label{thm:Sylows-completeness}
Let $G$ be a group satisfying the hypothesis of Problem~\ref{Prob:Lucc-type}, where $\mathcal{P}$ is the property of being complete. Then $G$ is metacyclic. 
\end{proposition}

\begin{proof}
If $G$ is non-metacyclic, then by Lemma \ref{lem:metacyclic}, $G$ contains a subgroup isomorphic to either $C_p \times C_p$ for some prime $p$ or the quaternion group $Q_8$. However, this contradicts the hypothesis in view of Lemmas \ref{lem:properties} (1) and Lemma \ref{lema:comprime}. The proof is complete. 
\end{proof}

We end this section with the following problem.
\begin{problem}
    Can we characterize all finite groups $G$ for which $\mathscr{P}_\bullet(G)$ is complete?
\end{problem}

\section{Connectivity}

Let  $G = \langle a \rangle \times \langle b \rangle$ be the elementary abelian $p$-group for some prime $p$. It is straightforward to verify that the element $(a,b)$ is an isolated vertex in the graph ${\mathscr{P}}^{*}_{\bullet}(G)$. Hence, this graph is disconnected.  However, the following result shows that every proper generalized power graph of a finite group whose center is not a $p$-group is connected.

\begin{theorem}\label{thm:Z(G)not-pgrooup}
If $G$ is a finite group such that $Z(G)$ is not a $p$-group, then  ${\mathscr{P}}^{*}_{\bullet}(G)$  is connected  with diameter  at most 4.
\end{theorem}
\begin{proof}
First, assume that $G$ is a finite abelian group. Since $\mathscr{P}^*_\bullet(G)$ coincides with the proper power graph of $G$ when $G$ is abelian, the connectedness statement follows from \cite[Theorem 2.4]{DFG}. For the reader's convenience, we provide a detailed proof below.

 Let $p$ be a prime dividing the order of $G$. By the Fundamental Theorem of Finite Abelian Groups, we can write $G = P \times Q$, where $P$ is a $p$-group and $Q$ is a non-trivial $p'$-group of order $q$.

Consider $(x,y), (x_1,y_1) \in G$ satisfying 
$$
|x| = p^n\neq 1, \quad |y|\neq 1, \quad |x_1| = p^{n_1}\neq 1, \quad |y_1|\neq 1.
$$
By Bezout's lemma, there exist integers $r$ and $s$ such that $r p^n + s q = 1$. Consequently,
\begin{align*}
(x,y)^{r p^n} &= (x^{r p^n}, y^{r p^n}) = (1, y^{r p^n})= (1, y^{1-sq})= (1, y), \\
(x,y)^{s q} &= (x^{s q}, y^{s q}) = (x^{s q}, 1)= (x^{1-rp^n}, 1)=(x,1).
\end{align*}

Thus, the vertex $(x,y)$ is connected to $(x,1)$ and $(1,y)$. Similarly, $(x_1,y_1)$ is connected to $(x_1,1)$ and $(1,y_1)$, and $(x_1,y)$ is connected to $(x_1,1)$ and $(1,y)$. Hence, we have the path $$(x,y) \to (1,y) \gets (x_1,y) \to (x_1,1) \gets (x_1,y_1),$$ which yields a path of length four between $(x,y)$ and $(x_1,y_1)$, which concludes the proof for the abelian case.

Now, assume that $G$ is a finite nonabelian group. 
Let $x$ and $y$ be non-trivial elements of $G$. 
If $[x,y] \neq 1$, then by Lemma \ref{L1}(8) we have $x \to [x,y] \leftarrow y$, which yields a path of length two between $x$ and $y$.
Thus, we may assume that $x$ and $y$ commute. Consider the subgroup $H = \langle x, y, Z(G) \rangle$, which is abelian. By assumption, $H$ is not a $p$-group, implying that $\mathscr{P}^*_\bullet(H)$ is connected with diameter at most $4$. This completes the proof.
\end{proof}

Next, we establish the following result.

\begin{theorem}\label{thm:centerless}
\begin{enumerate}
    \item If $G$ is a torsion group with trivial center, then
${\mathscr{P}}^{*}_{\bullet}(G)$ is connected with diameter at most $4.$
\item If $G$ is a finite group with a cyclic center, then $\mathscr{P}^*_\bullet(G)$ is connected with diameter at most $4$.
\end{enumerate}
\end{theorem}
\begin{proof}
\begin{enumerate}
    \item Let $x$ and $ y $ non-trivial elements of $G$. If $x$ and $y$ do not commute, then by Lemma \ref{L1}(8) we get $x\rightarrow [x,y]\leftarrow y$. Suppose now that $x$ and $y$ commute. Since neither element is central, there exist elements $x_1, y_1 \in G$ such that $x$ does not commute with $x_1$ and $y$ does not commute with $y_1$. By Lemma \ref{L1}(8) again, if either $x$ does not commute with $y_1$ or $y$ does not commute with $x_1$, then $x$ and $y$ are connected by a path of length at most $4$. Otherwise, $x$ commutes with $y_1$ and $y$ commutes with $x_1$. In this case, applying commutator calculus (cf. \cite[5.1.5 (ii)]{Rob}), we obtain the following path
    $$y\rightarrow [y,x_1y_1]\leftarrow x_1y_1\rightarrow [x,x_1y_1]\leftarrow  x.$$
     
    \item If $Z(G)$ is not a $p$-group, then Theorem \ref{thm:Z(G)not-pgrooup} applies and the proof is complete. Hence, we may assume that $Z(G)$ is a $p$-group. 
    Repeating the above argument, we obtain a path of length at most $5$ between any two non-central elements. Furthermore, since the center of $G$ is cyclic, say $Z(G)=\langle z \rangle$, any two elements of $Z(G)$ are connected by a path of length at most $2$. Suppose that there exists an element $g \in G$ whose order is coprime to $p$. 
The subgroup $\langle g, Z(G) \rangle$ is not a $p$-group. By Theorem \ref{thm:Z(G)not-pgrooup}, there exists a path of length at most 4 connecting $g$ to any element of $Z(G)$. Therefore, it remains to consider the situation where $G$ itself is a $p$-group.
In this case, let $x,y\in G\setminus \{1\}$ be such that $x$ is non-central and $y$ is central, and set $N = \langle x^G \rangle$, the normal closure of $x$. 
If $N = G$, then by Lemma~\ref{L1} (5) every element of $G$ is adjacent to $x$ in ${\mathscr{P}}^{*}_{\bullet}(G)$. 
If $N \neq G$, since $G$ is nilpotent, we have $N \cap Z(G) \neq \{1\}$. Choosing a non-trivial element $z^i$ in this intersection, we obtain the following path connecting $x$ and $y$ $$x\rightarrow z^i\leftarrow z \rightarrow y.$$ This completes the proof.\end{enumerate} 
\end{proof}

Note that the previous theorem ensures that the proper generalized power graph of groups such as the Gupta-Sidki, Grigorchuk, and Tarski monster groups is connected.

\begin{problem}
  Can the ``finite group'' hypothesis be replaced with ``torsion group'' in Item (2) of Theorem \ref{thm:centerless}?
\end{problem}

In view of the two preceding theorems, it remains to consider the case where the center is a non-cyclic $p$-group. To this end, we present a finite group $G$ with such a center for which $\mathscr{P}^*_\bullet(G)$ is disconnected.

\begin{example}
Let $G = D_8 \times C_2$ be the direct product of the dihedral group $D_8 = \langle a, b \mid a^4 = 1, b^2 = 1, bab = a^{-1} \rangle$ and the cyclic group $C_2 = \langle c \mid c^2 = 1 \rangle$. 
We claim that $(1,c)$ is an isolated vertex of $\mathscr{P}^*_\bullet(G)$. Indeed, the center of $G$, given by $Z(G) = \{ (1, 1), (a^2, 1), (1, c), (a^2, c) \}$, is isomorphic to the Klein four-group.  
It is straightforward to observe that if $(x,y)$ is adjacent to $(1,c)$, then $(x,y) \to (1,c)$. In this case, there exist elements $(g_1, h_1), \ldots, (g_n, h_n) \in G$ such that$$(1,c) = \prod_{i=1}^{n} (x,y)^{(g_i,h_i)} = \left( \prod_{i=1}^{n} x^{g_i}, y^n \right).$$ 
Consequently, $y = c$, $n$ is odd, and $\prod_{i=1}^{n} x^{g_i} = 1$. Since $\vert{}D_8\vert{} = 8$, the element $x$ must be non-central. Thus, $x$ belongs to one of the following conjugacy classes:
$$b^{D_8} = \{b, a^2b\}, \quad (ab)^{D_8} = \{ab, a^3b\}, \quad \text{and} \quad a^{D_8} = \{a, a^3\}.$$However, any two elements within the same conjugacy class commute. Thus, no such element $x$ exists, which implies that  $\mathscr{P}^*_\bullet(G)$ is disconnected. Direct calculations show that the graph $\mathscr{P}^*_\bullet(G)$ can be represented as follows:

\end{example}

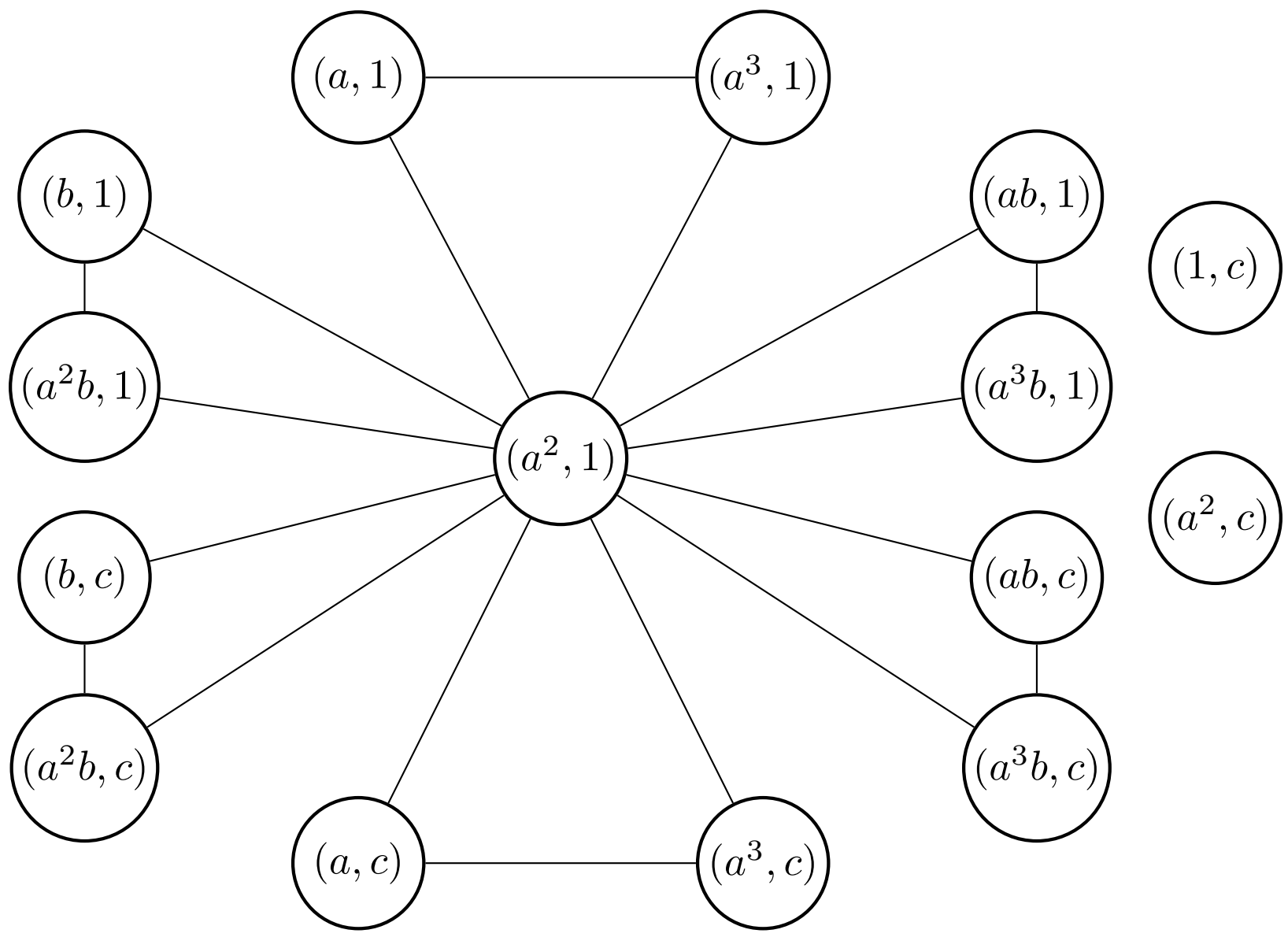
\begin{center}
\begin{tikzpicture}[
  vertex/.style={circle, draw, thick, minimum size=1.1cm, inner sep=1pt},
  every node/.style={font=\small}
]

\node[vertex] (A) at (0,0) {$(a^2,1)$};

\node[vertex] (B1) at (-1.7,3.2) {$(a,1)$};
\node[vertex] (B2) at ( 1.7,3.2) {$(a^3,1)$};

\node[vertex] (C1) at (-4.0,2.2) {$(b,1)$};
\node[vertex] (C2) at (-4.0,0.6) {$(a^2b,1)$};

\node[vertex] (D1) at (4.0,2.2) {$(ab,1)$};
\node[vertex] (D2) at (4.0,0.6) {$(a^3b,1)$};

\node[vertex] (E1) at (-4.0,-1.0) {$(b,c)$};
\node[vertex] (E2) at (-4.0,-2.6) {$(a^2b,c)$};

\node[vertex] (F1) at (4.0,-1.0) {$(ab,c)$};
\node[vertex] (F2) at (4.0,-2.6) {$(a^3b,c)$};

\node[vertex] (G1) at (-1.7,-3.4) {$(a,c)$};
\node[vertex] (G2) at ( 1.7,-3.4) {$(a^3,c)$};

\node[vertex] (H1) at (5.5,1.6) {$(1,c)$};
\node[vertex] (H2) at (5.5,-0.5) {$(a^2,c)$};

\draw (B1)--(B2);
\draw (C1)--(C2);
\draw (D1)--(D2);
\draw (E1)--(E2);
\draw (F1)--(F2);
\draw (G1)--(G2);

\draw (A)--(B1);
\draw (A)--(B2);
\draw (A)--(C1);
\draw (A)--(C2);
\draw (A)--(D1);
\draw (A)--(D2);
\draw (A)--(E1);
\draw (A)--(E2);
\draw (A)--(F1);
\draw (A)--(F2);
\draw (A)--(G1);
\draw (A)--(G2);

\end{tikzpicture}
\end{center}

\begin{theorem}
    Let $n\geq 2$ be an integer and $G=SL_n(F)$ be the special linear group of degree $n$ over a field $F$. Then ${\mathscr{P}}^{\ast}_{\bullet}(G)$  is connected. 
\end{theorem}
\begin{proof}
Assume that $n>2$ or $\vert{}F\vert{}>3$, and let $A\in G$ be a non-scalar matrix. By \cite{Wit}, every element $B$ in $G$ is a product of conjugates of $A$, and so $A\rightarrow B$. Thus, in order to complete the proof, it is sufficient to show that ${\mathscr{P}}^{\ast}_{\bullet}(G)$ is connected when $G=SL(2,2)$ and $G=SL(2,3)$.The first case, $G=SL(2,2)$, follows from Theorem \ref{thm:centerless} (1) because $SL(2,2)\simeq S_3$. Now, we will deal with the case $G=SL(2,3)$. Here, $G$ has order $24$ and contains a normal subgroup isomorphic to the quaternion group $Q_8$. By Theorem \ref{thm:centerless} (2) (or the figure in Remark \ref{remark:order8}), ${\mathscr{P}}^{\ast}_{\bullet}(Q_8)$ is connected. 
Further, $G$ has exactly four Sylow $3$-subgroups, generated respectively by the matrices$$A=\begin{pmatrix} 1 & 1 \\ 0 & 1 \end{pmatrix}, \quad B=\begin{pmatrix} 1 & 0 \\ 1 & 1 \end{pmatrix}, \quad C=\begin{pmatrix} 0 & 2 \\ 1 & 2 \end{pmatrix}, \quad D=\begin{pmatrix} 2 & 2 \\ 1 & 0 \end{pmatrix}.$$
Since all Sylow 3-subgroups are conjugate, all elements of order 3 are connected. Consider arbitrarily $M$ a non-trivial element in $G$. Note that  $G^3\leq Q_8$. If $M^3\neq 1$, then $M$ is connected to all elements in $Q_8$. If $M^3=1$, then $M$ is contained in some 3-Sylow subgroup. Further, the commutator $[A,B] = \begin{pmatrix} 0 & 1 \\ -1 & 1 \end{pmatrix}$ has order $4$, and thus $[A,B]$ is connected to all elements in $Q_8$. By Lemma \ref{L1} (8), $C \rightarrow A \rightarrow [A,B] \leftarrow B \leftarrow D$. This shows that every non-trivial element of $G$ is connected to all non-trivial elements of $Q_8$. This completes the proof.
\end{proof}

\begin{theorem}\label{thm:Sylows-connectivity}
\begin{enumerate}
    \item Let $G$ be a finite group of odd order satisfying the hypothesis of Problem~\ref{Prob:Lucc-type}, where $\mathcal{P}$ is the property of being connected. Then $G$ is metacyclic.
    \item If $G$ is a finite group such that ${\mathscr{P}}^*_{\bullet}(H/N)$ is connected for every proper subgroup $H$ of $G$ and every  $N\triangleleft H$, then  $G$ is metacyclic.
\end{enumerate}
\end{theorem}
\begin{proof}
\begin{enumerate}
    \item     Assume that $G$ is a non-metacyclic group. Since $G$ has odd order, Lemma~\ref{lem:metacyclic} implies that $G$ contains a subgroup isomorphic to $C_p \times C_p$ for some prime $p$. However, because $\mathscr{P}^*_{\bullet}(C_p \times C_p)$ is not connected, this yields a contradiction.

\item Suppose that $G$ is a non-metacyclic group. Lemma~\ref{lem:metacyclic} implies that $G$ contains a subgroup isomorphic to $C_p \times C_p$ for some prime $p$ or to a generalized quaternion group. 
First, assume that $G$ is a $p$-group for some prime $p$. If $p=2$, our assumption implies that $G$ has a unique subgroup of order 2. By \cite[5.3.6]{Rob}, $G$ is either cyclic or a generalized quaternion group, yielding a contradiction in both cases. If $p$ is odd, Lemma~\ref{lem:metacyclic} implies that $G$ contains a subgroup isomorphic to $C_p \times C_p$, which is again a contradiction. 

Now, assume that $G$ is not a $p$-group. By Item (1), we may assume that $G$ has even order. Furthermore, it follows from \cite[5.3.6]{Rob} that every Sylow subgroup of odd order in $G$ is cyclic, and every Sylow $2$-subgroup is either cyclic or a generalized quaternion group. 
Let $P$ be a Sylow $2$-subgroup of $G$. If $P$ is cyclic, then \cite[Theorem 7.6.2]{Gor} implies that $G$ is metacyclic, a contradiction. Consequently, $P$ is a generalized quaternion group. But $P/P' \simeq C_2 \times C_2$, which is a contradiction. This completes the proof.
\end{enumerate}
\end{proof}

\section{Planarity}
In this section, we provide necessary and sufficient conditions for the planarity of the proper generalized  power graph ${\mathscr{P}}^{*}_{\bullet}(G)$ of a finite group $G$. More precisely, we prove the following result.

\begin{theorem}\label{thm:planar}
Let $G$ be a non-trivial finite group. Then ${\mathscr{P}}^{*}_{\bullet}(G)$ is planar if and only if one of the following conditions holds:
    \begin{enumerate}
    \item $G\simeq C_3^m$, for some integer $m\geq 1$;
         \item $G\simeq C_5^m$, for some integer $m\geq 1$;
                 \item $G\simeq C_2^m\times C_3^n$, where either $m=1$ or $n=1$;
               \item $G\simeq  C_2^m\times C_4^n$, for some non-negative integers $m$ and $n$, not both zero;
               \item $G$ is a non-abelian group of exponent $4$ that is a central extension of $C_2$ by $C_2^m$ for some positive integer $m$.
        \end{enumerate}
\end{theorem}

To prove the Theorem \ref{thm:planar}, we first establish several auxiliary lemmas.

\begin{lemma}\label{lem:C_6}
The graph ${\mathscr{P}}^{*}_{\bullet}(C_6\times  C_6)$ is non-planar.
\end{lemma}
\begin{proof}
Consider the additive group $C_6 \times C_6$. Let $a_1=(3,0)$, $a_2=(0,3)$, and $a_3=(3,3)$ be elements of order 2, and let $b_1=(0,2)$, $b_2=(0,4)$, and $b_3=(2,0)$ be elements of order 3. For each $i, j \in \{1,2,3\}$, define $c_{ij} = a_i + b_j$.
Note that $4c_{ij} = b_j$ and $3c_{ij} = a_i$, which implies that for each pair $(i,j)$ there exists a path $b_j \leftarrow c_{ij} \to a_i$. Consequently, the graph $\mathscr{P}^*_\bullet(C_6 \times  C_6)$ contains a subdivision of $K_{3,3}$. By Kuratowski's theorem \cite[Theorem 9.10, p. 153]{Bon}, $\mathscr{P}^*_\bullet(C_6 \times C_6)$ is non-planar, which completes the proof.
\end{proof}

Let $\pi$ be a set of prime numbers. Recall that a group $G$ is called a $\pi$-group if  every non-trivial element in $G$ has an order whose prime factors belong to $\pi$.

\begin{lemma}\label{obs}
Let $G$ be a non-trivial finite group such that ${\mathscr{P}}^{*}_{\bullet}(G)$ is planar. Then $G$ is a metabelian $\{2,3,5\}$-group and every conjugacy class of $G$ has size at most 4.
\end{lemma}

\begin{proof}
    Suppose there exists a prime $p \geq 7$ dividing the order of $G$. 
Then, by Cauchy's Theorem, $G$ contains a cyclic subgroup $C$ of order $p$. 
Consequently, 
${\mathscr{P}}^{*}_{\bullet}(C) \simeq K_{p-1}$.
Since $p-1 \geq 6$, it follows that ${\mathscr{P}}^{*}_{\bullet}(G)$ contains $K_5$ as a subgraph. 
By Kuratowski's Theorem, ${\mathscr{P}}^{*}_{\bullet}(G)$ is non-planar, which implies that $G$ must be a $\{2,3,5\}$-group.

Now let $x$ be a non-trivial element of $G$. If $|x^G| \geq 5$, then by Lemma \ref{L1} (3), 
${\mathscr{P}}^{*}_{\bullet}(G)$  again contains $K_5$ as a subgraph, contradicting planarity. 
Hence, we have that every conjugacy class of $G$ has size at most $4$. 
Under this assumption, it follows from \cite[Theorem 5.6 and Section 6]{W} that the derived subgroup $G'$ is abelian, completing the proof.
\end{proof}

\begin{lemma} \label{Caso3}
Let $G$ be a finite group containing an element $a$ with at least three conjugates. If ${\mathscr{P}}^{*}_{\bullet}(G)$ is planar, then $\langle a^G \rangle = a^G \cup \{1\}$. Furthermore, $\langle a^G \rangle$ is isomorphic to $V_4$, the Klein four-group.
\end{lemma}
\begin{proof}
First, suppose that $a$ is adjacent to some element $b \neq 1$ with $b \notin a^G$. We claim that $b$ is a central element of order $2$. Indeed, if $b \notin Z(G)$, then $\vert{}b^G\vert{} \ge 2$. By Lemma \ref{L1} (3) and (4), every element in $a^G$ is adjacent to every element in $b^G$ and thus $K_5$ is a subgraph of ${\mathscr{P}}^{*}_{\bullet}(G)$, which contradicts Kuratowski's theorem.  Therefore, $b \in Z(G)$. Furthermore, if $b^2 \neq 1$, then every element in $a^G$ is adjacent to every element in $\langle b \rangle$, which again yields that $K_5$ is a subgraph of ${\mathscr{P}}^{*}_{\bullet}(G)$, reaching a contradiction. 

Now, suppose that there exist distinct non-trivial elements $b_1, b_2, \dots, b_n \notin a^G$, with $n \ge 3$, adjacent to $a$. Lemma \ref{L1}(3) and (4) imply that $K_{3,3}$ is a subgraph of ${\mathscr{P}}^{*}_{\bullet}(G)$, again contradicting planarity.
Thus, the order of $\langle a^G \rangle$ can only be 4, 5, or 6. If $|\langle a^G \rangle| = 6$, then there exist two central involutions $b_1, b_2 \in \langle a^G \rangle$ such that $b_1, b_2 \notin a^G$. Since 3 divides $|\langle a^G \rangle|$, the subgroup $\langle a^G \rangle$ contains an element of order 3; because $b_1$ and $b_2$ are involutions, $a$ must have order 3. However, no group of order 6 contains two central involutions (or three elements of order 3), so $|\langle a^G \rangle| \neq 6$. If we suppose that $|\langle a^G \rangle| = 5$, there exists a central involution $b \in \langle a^G \rangle$, which contradicts Lagrange's theorem. Thus, the order of $\langle a^G \rangle$ must be 4. Since all elements in $a^G$ are conjugate and therefore have the same order, it follows that $\langle a^G \rangle \cong V_4$, which completes the proof.
\end{proof}

\begin{lemma}\label{case4}
Let $G$ be a finite group in which no element has more than 4 conjugates, and some element has exactly 4 conjugates. If ${\mathscr{P}}_{\bullet}^*(G)$ is planar, then $G$ has some element of order 5. 
\end{lemma}
\begin{proof}
It is sufficient to show that $G$ contains a subgroup of order $5$. Let $x$ be an element with exactly $4$ conjugates, and consider $H = \langle x^G \rangle$, the normal closure of $x$. It is clear that $|H| \ge 5$. Suppose, for a contradiction, that $|H| \ge 6$. Then there exists a non-trivial element $y \in H \setminus x^G$. By Lemma~\ref{L1}(4), $y$ is adjacent to all elements of $x^G$, implying that ${\mathscr{P}}_{\bullet}(G)$ contains $K_5$ as a subgraph. By Kuratowski's Theorem \cite{Bon}, ${\mathscr{P}}_{\bullet}(G)$ is non-planar, contradicting our assumption. Thus, $|H| = 5$, completing the proof.
\end{proof}

\begin{lemma}\label{lemma-to-planarity}
Let $G$ be a non-trivial finite group.
\begin{enumerate}
   \item   If $G$ is an extension of $C_3$ by $C_5^m, m\geq 1$, then $G$ is abelian. Consequently, ${\mathscr{P}}^{*}_{\bullet}(G)$ is non-planar.
    \item   If $G$ is an extension of $N$ by $C_5^m, m\geq 1$, where $|N|=2^n, n\leq 3$, then ${\mathscr{P}}^{*}_{\bullet}(G)$ is non-planar.
\end{enumerate}
  
\end{lemma}
\begin{proof}
We only present the proof of (1), since (2) follows similarly. By Sylow's theorem, $G$ contains a unique Sylow 5-subgroup, which is therefore normal. Since its Sylow subgroups are normal and abelian with coprime orders, $G$ is an abelian group containing an element of order 15. The result then follows from \cite[Theorem 4.2]{DR}.
\end{proof}

\begin{proof}[{\bf Proof of Theorem \ref{thm:planar}}]
First, suppose that ${\mathscr{P}}^{*}_{\bullet}(G)$ is planar. 
By Lemma \ref{obs} we may assume henceforth that $G$ is a metabelian $\{2,3,5\}$-group and every conjugacy class of $G$ has size at most $4$. We prove this direction of the theorem by considering each case separately.

{\bf Case 1:} Suppose that every conjugacy class has cardinality 1. In this case, $G$ is abelian. If $G$ has a non-trivial Sylow $5$-subgroup, then by Theorem 4.2 of \cite{DR}, 
we have 
$G \simeq C_5^m$
for some positive integer $m$ (since the existence of any other Sylow subgroup would force $G$ to contain 
a cyclic subgroup of order greater than 6). Otherwise, 
$G = P_2 \times P_3$, where $P_2$ and $P_3$ are Sylow $p$-subgroups for $p = 2, 3$, respectively (each possibly trivial). 
 If $P_2 = 1$, then $G \simeq C_3^n$
for some $n \geq 1$. If $P_3 = 1$, then $G$ has exponent dividing $4$.
If both are non-trivial, Theorem 4.2 of \cite{DR} yields 
$G \simeq C_2^m \times C_3^n$ for positive integers $m, n$. By Lemma \ref{lem:C_6}, if $m\geq 2$ and $n\geq 2$, then ${\mathscr{P}}^{*}_{\bullet}(G)$ is non-planar. Thus, either $m=1$ or $n=1.$

{\bf Case 2:} Suposse that every conjugacy class has cardinality at most 2, and there exists an element with exactly two conjugates. 
By Theorem 5.6 of \cite{W}, the derived subgroup $G'$ is cyclic of order $2$, say $G'=\langle z\rangle$. Since $z^g\in G'=\{1,z\}$ for all $g$ in $G$, we get that 
$G' \leq  Z(G) $. Thus, we conclude that $G$ is nilpotent, and therefore 
it decomposes as the direct product of its Sylow subgroups.  By \cite[Theorem 4.2]{DR}, the Sylow 5-subgroup of $G$ is trivial. Furthermore, \cite[Theorem 5.6]{W} also implies that the Sylow 3-subgroup of $G$ is trivial.

Since $G$ is nonabelian, it follows from \cite[Theorem 4.2]{DR} that $G$ has exponent exactly 4. Thus, $G$ is a extension of $C_2$ by $ C_2^m \times C_4^n$, for some integers $m$ and $n$. Next, we show that $n = 0$. To this end, suppose by contradiction that $n \ge 1$.

Firstly, note that if $x$ is a non-central element of $G$, then $x^g = x[x,g]$ for all $g \in G$. Since $|x^G| = 2$ and $G'=\{1, z\}$, it follows that $x^G = \{x, xz\}$.

Now, set $T = \{g \in G \mid g^2 \notin G'\}$. Since $n \ge 1$, we have $T \neq \emptyset$. Note that $G \setminus T$ is a proper subgroup of $G$, which implies that $G = \langle T \rangle$. If $T \subseteq Z(G)$, then $G$ is abelian, a contradiction.
Thus there exist $y\in T\setminus Z(G)$. By previous paragraph, we have that $y^G=\{y,yz\}$ and  $$\langle y^G\rangle=\langle y,z\rangle=\{1,z,y,y^2,y^3,yz,y^2z,y^3z\}.$$
Note that since $y$ is noncentral, $y^3=y^{-1}$, $yz$ and $y^3z$ are also noncentral.

Thus, by the same argument as above, we have that  $(yz)^G=\{yz,y\}$, \linebreak ${(y^3)}^{G}=\{y^3,y^3z\}={(y^3z)}^G$. Therefore, we conclude that $$\langle y^G\rangle=\langle (yz)^G\rangle=\langle (y^3)^G\rangle=\langle (y^3z)^G\rangle=\langle y,z\rangle.$$ 
Consequently, $K_5$ is a subgraph of ${\mathscr{P}}^{*}_{\bullet}(G)$, a contradiction. Thus $n = 0$, which implies that $G$ has exponent 4 and is an extension of $C_2$ by $C_2^m$ for some integer $m$.

{\bf Case 3:}
Suppose that every conjugacy class in $G$ has cardinality at most 3, and there exists an element with exactly 3 conjugates.
By \cite[Theorem 5.6]{W}, the derived subgroup $G'$ is cyclic of order $3$.

Let  $S = \{x \in G \mid |x^G| = 3\}$ be the set of all elements in $G$ whose conjugacy classes have size $3$. For each $x \in S$, set $H_x = \langle x^G \rangle$, and define the subgroup $H =\prod_{x\in S} H_x $. Consider the quotient group $G/H$. Note that every conjugacy class of $G/H$ has size at most $2$. Further, its derived subgroup $(G/H)' \cong {G'}/{H \cap G'}$ has order dividing 3.

Thereby, if there exists an element in $G/H$ with exactly 2 conjugates, then \cite[Theorem 5.6]{W} ensures that $(G/H)'$ is cyclic of order 2, a contradiction. Hence, $G/H$ must be abelian, so $(G/H)' = \{1\}$, which implies $G' \subseteq H$. However, by Lemma \ref{Caso3}, $H$ is a 2-group, yielding a final contradiction.

{\bf Case 4:} Suppose that every conjugacy class in $G$ has cardinality at most $4$ and there exists an element with exactly 4 conjugates.
According to Section 6 of \cite{W}, the derived subgroup $G'$ is isomorphic to $C_4$, $C_2^2$, or $C_2^3$. Moreover, Lemma \ref{case4} implies that 5 divides the order of $G$. Hence, by Theorem 4.2 of \cite{DR}, we obtain $G/G'\simeq C_5^m$ for some integer $m\geq 1$. However, Lemma \ref{lemma-to-planarity} (2) guarantees that $G$ is abelian and contains a cyclic subgroup of order $10$, which contradicts Theorem 4.2 of \cite{DR}. 

Now, we provide a proof of the converse. If $G$ satisfies one of the first three conditions, then \cite[Theorems 4.3 and 4.5]{DR} ensure that ${\mathscr{P}}^{*}_{\bullet}(G)$ is planar. Theorem 2 of \cite{AAA} implies that ${\mathscr{P}}^{*}_{\bullet}(G)$ is planar if $G$ satisfies the conditions of item (4). Finally, assume that $G$ satisfies the conditions of item (5). 
Thus, there exists a central involution $z$ satisfying $G/\langle z\rangle$ is elementary abelian 2-group. Since $G$ is nonabelian, we get $G'=\langle z\rangle$. Consequently, for any pair of elements $x,g\in G$ we get $x^g=x[x,g]$ which implies that $x^G=\{x,xz\}$.

Assume that $y \rightarrow x$. We show that $x \in \{z, yz\}$. Indeed, there exist elements $g_1, \ldots, g_n \in G$ satisfying $$x = \prod_{i=1}^n y^{g_i} \in \langle y^G \rangle = \langle y, yz \rangle = \langle y, z \rangle.$$ 
Since $z$ is a central involution, we have either $x = y^j$ or $x = y^j z$ for some integer $j$. Furthermore, since $y^2 \in \langle z \rangle = \{1, z\}$, it is straightforward to see that in both cases $x \in \{z, yz\}$. This shows that if $y$ is not an isolated vertex, then its unique adjacent vertices are $z$ and $yz$. Up to isolated vertices, the graph $\mathscr{P}^*_\bullet(G)$ can be represented as shown in the figure below. 
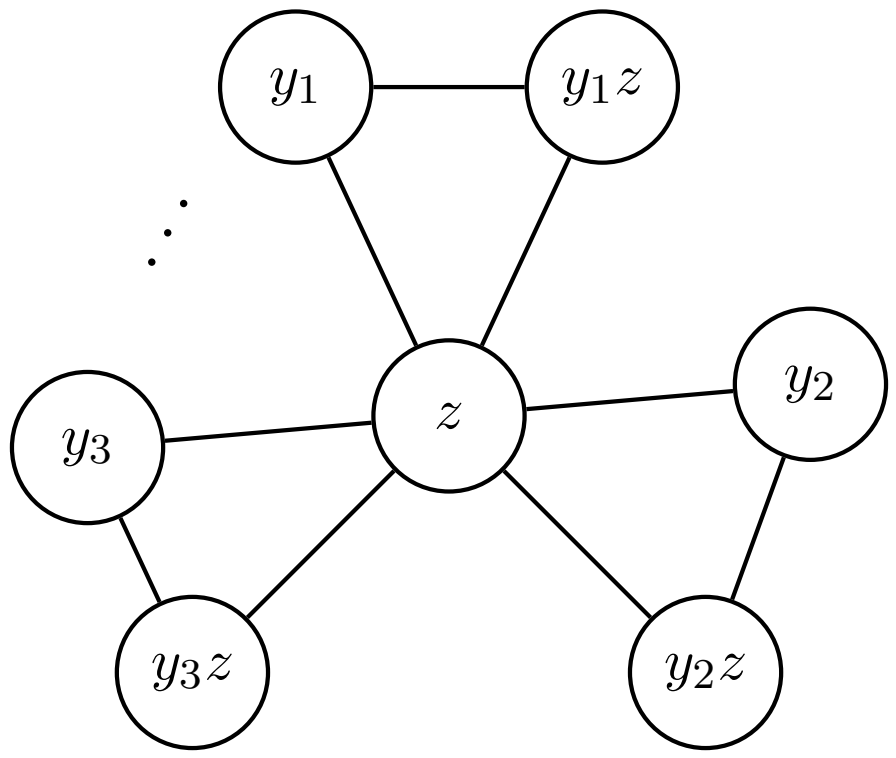
\begin{figure}[H]
\centering
\begin{tikzpicture}[
    vertice/.style={
        circle, 
        draw, 
        thick, 
        minimum size=10mm, 
        inner sep=0pt, 
        font=\normalsize,
        fill=white
    },
    aresta/.style={draw, thick}
]

    \node[vertice] (z) at (0,0) {$z$};

    \node[vertice] (y1)  at (115:2.4) {$y_1$};
    \node[vertice] (y1z) at (65:2.4)  {$y_1z$};

    \node[vertice] (y2)  at (5:2.4)   {$y_2$};
    \node[vertice] (y2z) at (-45:2.4) {$y_2z$};

    \node[vertice] (y3z) at (-135:2.4) {$y_3z$};
    \node[vertice] (y3)  at (-175:2.4) {$y_3$};

    \node[font=\Large, rotate=90] at (145:2.4) {$\ddots$ \ \ };

    \draw[aresta] (z) -- (y1);
    \draw[aresta] (z) -- (y1z);
    \draw[aresta] (y1) -- (y1z);

    \draw[aresta] (z) -- (y2);
    \draw[aresta] (z) -- (y2z);
    \draw[aresta] (y2) -- (y2z);

    \draw[aresta] (z) -- (y3);
    \draw[aresta] (z) -- (y3z);
    \draw[aresta] (y3) -- (y3z);

\end{tikzpicture}

\end{figure}
This completes the proof.

\end{proof}

\section*{IA}
Artificial intelligence (https://gemini.google.com/app) was used to assist with routine typesetting, figure generation, grammar checking, and phrasing. Apart from these minor  editorial and illustrative uses, the text and results of this paper were entirely authored by the human authors.

\section*{Disclosure statement}
The authors have no conflicts of interest to declare that are relevant to the content of this article.

\section*{Data sharing}
Not applicable to this article as no datasets were generated or analyzed during the current study.

\section*{Funding}
This work was partially supported by FAPEMIG (grant number RED-00133-21).

\section*{Acknowledgement}
The second author acknowledges the Departamento de Matemática of the Universidade Federal de Minas Gerais for its excellent hospitality while this research was conducted.

\end{document}